\documentclass[12pt,a4paper,leqno]{article}

\usepackage{amsmath}
\usepackage{amsfonts}

\newtheorem{theorem}{Theorem}
\newtheorem{definition}[theorem]{Definition}
\newtheorem{lemma}[theorem]{Lemma}
\newtheorem{remark}[theorem]{Remark}

\newenvironment{proof}[1][Proof]{\noindent\textbf{#1.} }{\ \rule{0.5em}{0.5em}}

\makeatletter
\@addtoreset{equation}{section}

\makeatother 
\newdimen\dummy
\dummy=\oddsidemargin
\begin{document}

\title{The geometry of contact metric three manifolds}
\author{Karatsobanis J.\\
Mathematics Division\\
School of Technology\\
Aristotle University of Thessaloniki\\
Thessaloniki 54124 - Greece\\
glenn@gen.auth.gr\\
Math. subj. class [2000] 53D10\\
}

\maketitle


\begin{abstract}
Determining the associated metrics we get a local classification of contact metric three manifolds.
\end{abstract}

\section{Introduction}\label{intro}

It is known that a $C^{\infty}$ manifold $M$ admits many Riemannian metrics. Once such a metric has been specified, $M$ has been given a great rigidity. This rigidity imposed by the metric on $M$ distinguishes the subject of geometry and topology.

Riemannian metrics associated to a contact or symplectic structure have been studied extensively. They form an infinite dimensional space and each element of this space defines the same volume element.

The associated metrics of a three dimensional contact metric manifold are determined in the present paper through an explicit formula ((\ref{metric-simp}), p. \pageref{metric-simp}). This formula comes from the solution of a partial differential equation system ((\ref{S1}) - (\ref{S5}), p. \pageref{S1}) which we call \emph{contact system} (Definition \ref{cs1}, p. \pageref{cs1}). We introduce a special coordinate system (simplifying coordinates,
Definition \ref{simpdef}, p. \pageref{simpdef}) such that the contact system is reduced to a
Riccati equation (\ref{S*3}, p. \pageref{S*3}).

Every contact manifold carries a contact metric structure, i.e. a Riemannian metric compatible with the contact structure \cite{B02}. It is a classic result \cite{MA71,LU70} that every $3$-dimensional closed and orientable manifold admits a contact structure and hence, a contact metric structure. Thus the metric (\ref{metric-simp}) represents the local geometry of the contact metric structure of these manifolds.

Concerning the existence of contact structures on open manifolds the result of Gromov \cite{G86} states that a non-compact, odd-dimensional manifold carrying a
hyperplane field with an almost complex structure, accepts a contact structure. Here we provide an existence theorem (section \ref{example}, Theorem \ref{lemma1}) for contact metric structures on $3$-dimensional orientable manifolds, both compact and non-compact. The proof of the Theorem \ref{lemma1} founds the above mentioned contact system (Definition \ref{cs1}) whose solution leads to the determination of the associated metrics. This theorem requires a Levi-Civita connection with a particular framing. Its premises are indeed met since there are many explicit examples of contact metric $3$-manifolds satisfying them, e.g. the example at section \ref{example1} of this paper.

In section \ref{algorithm} we give an algorithm for the construction of contact metric structures in $3$-manifolds. Finally in section \ref{example1} we give an example of a contact metric $3$-manifold utilizing the algorithm and the simplifying coordinates.

\section{Preliminaries}\label{intro2}

In what follows if $f\left( x^{1},x^{2},\ldots ,x^{n}\right) $ is a smooth function then $f_{i}$ will denote the partial derivative $\frac{\partial f}{\partial x^{i}}$. If $B$ is some $n\times n$ invertible matrix with elements $b_{\alpha }^{i}$, $\alpha =1,\ldots n$, $i=1,\ldots n$, then $\left(b^{-1}\right) _{i}^{\alpha }$, $\alpha =1,\ldots n$, $i=1,\ldots n$, will denote the element of the $i$ row and the $\alpha $ column of the matrix $B^{-1}$.

In this section we collect some basic facts on contact metric manifolds. All manifolds are assumed to be smooth.

A differentiable $(2n+1)$-dimensional manifold $M$ is called a contact manifold if it carries a nowhere integrable affine distribution of hyperplanes. It can be shown that every contact manifold carries a $1$-form $\eta $ satisfying $\eta \wedge (d\eta )^{n}\neq 0$ everywhere on $M$. It is well known \cite{B02} that given $\eta$ there exists a unique vector field $\xi$ (called the characteristic vector field or Reeb field) defined by $(d\eta )(\xi ,X)=0$ and $\eta (\xi)=1$. Polarizing $d\eta $ on the contact $2n$-subbundle $D$ given by $\eta =0$, one obtains a Riemannian metric $g$ and a $(1,1)$-tensor field $\phi $ such that:
\begin{equation*}
(d\eta )(X,Y)=g(\phi X,Y),~~~\eta (X)=g(X,\xi ),~~~\phi ^{2}=-I+\eta \otimes
\xi 
\end{equation*}
The metric $g$ is called \emph{associated metric} of $\eta $, and $(\phi,\eta ,\xi ,g)$ a contact metric structure \cite{B02}. The associated metric is not unique. A differentiable $(2n+1)$-dimensional manifold $M$ equipped with a contact metric structure is called \emph{contact metric} $(2n+1)$- manifold.

In the theory of contact metric manifolds the symmetric tensor field $h:=\frac{1}{2}\mathcal{L}_{\xi }\phi $, $\mathcal{L}$ being the Lie derivative, plays a fundamental role. It satisfies
\begin{eqnarray}
h\xi  &=&0,~~~Tr\left( h\right) =0,~~~Tr(h\phi )=0,~~~h\phi =-\phi h,~~~
\label{1} \\
\eta \circ h &=&0,\text{ }(\nabla _{\xi }h)\phi =-\phi (\nabla _{\xi
}h),~~~(\nabla _{\xi }h)\xi =0  \notag
\end{eqnarray}
The equation $h=0$ holds if and only if $\xi $ is Killing with respect to the associated metric, and then $M$ is called $K$-\emph{contact}. On a contact metric manifold we have $h\phi X=\nabla _{X}\xi +\phi X$ and\ $\nabla _{\xi }\phi =0$, where $\nabla $ is the Levi-Civita connection and $X$ is any vector field.

If the almost complex structure $J$ on $M\times \mathbb{R}$ defined by $J(X,f\frac{d}{dt})=(\phi X-f\xi ,\eta (X)\frac{d}{dt})$ is integrable, then $M$ is said to be Sasakian. A Sasakian manifold is $K$-contact and the converse is true only for $3$-dimensional manifolds.

The sectional curvature of a plane section containing $\xi $ is called $\xi $-sectional curvature. The element of the contact subbundle $D$ over a point $p$ is denoted by $D_{p}$. If $X\in D_{p}$ we denote the $\xi $-sectional curvature by $K(X,\xi )$. The sectional curvature $K(X,\phi X)$ of a plain section spanned by the vector fields $X$ and $\phi X$ (both orthogonal to $\xi $) is called $\phi $-\emph{sectional curvature}.

If $e\in D_{p}$ is a unit vector then $\phi e$ is another unit vector of $D_{p}$ orthogonal to $e$ and $\left\{ e,\phi e\right\} $ is the so-called symplectic basis of $D_{p}$. The triple $\left\{ e,\phi e,\xi \right\} $ is a basis of $T_{p}M^{3}$ and is called $\phi $-basis. Since there is a rotational freedom in choosing the symplectic basis it follows
that if $h\neq 0$ at $p$, this basis can be chosen to diagonalize $h$. Such a $\phi $-basis is called $\phi $-eigenbasis. From (\ref{1}) it follows that if $he=\lambda e$ then $h\phi e=-\lambda \phi e$.

Let $M^{3}:=(M^{3},\phi ,\eta ,\xi ,g)$ be a $3$-dimensional contact metric manifold (\emph{contact metric $3$-manifold}). Let $U_{0}$ and $U_{1}$ be two open sets defined by:
\begin{equation*}
U_{0}:=\left\{ p\in M^{3}:\lambda =0\text{ in a neighborhood of }p\right\} 
\end{equation*}
\begin{equation*}
U_{1}:=\left\{ p\in M^{3}:\lambda \neq 0\text{ in a neighborhood of }
p\right\} 
\end{equation*}
The closure of $U_{0}\cup U_{1}$ equals to $M^{3}$, i.e. $U_{1}\cup U_{0}$ is open and dense in $M^{3}$. For every point $p\in U_{0}\cup U_{1}$ there exists a $\phi $-eigenbasis $\{e,\phi e,\xi \}$ in a neighborhood of $p$.

Using (\ref{1}) and assuming that $\left\{ e,\phi e,\xi \right\} $ is a $\phi $-eigenbasis, it is not hard to prove the following:

\begin{lemma}
\label{aux1}On a contact metric $3$-manifold $M^{3}$ we have:
\begin{align*}
\nabla _{e}e& =b\,\phi e,~~~\nabla _{e}\phi e=-b\,e+(\lambda +1)\,\xi
,~~~\nabla _{e}\xi =-(\lambda +1)\,\phi e \\
\nabla _{\phi e}e& =-c\,\phi e+(\lambda -1)\,\xi ,~~~\nabla _{\phi e}\phi
e=c\,e,~~~\nabla _{\phi e}\xi =(1-\lambda )\,e \\
\nabla _{\xi }e& =a\,\phi e,~~~\nabla _{\xi }\phi e=-a\,e
\end{align*}
\end{lemma}

\begin{remark}
Recalling that on $U_{0}$ the relation $\nabla _{X}\xi +\phi X=0$ holds, Lemma \ref{aux1} is true on $U_{0}\cup U_{1}$. But the relations in the Lemma \ref{aux1} are equations between continuous functions. Moreover $U_{0}\cup U_{1}$ is open and dense in $M^{3}$. Hence the Lemma \ref{aux1} can be extended to $M^{3}$.
\end{remark}

From Lemma \ref{aux1} and the well known property $\nabla _{X}Y-\nabla _{Y}X=\left[ X,Y\right]$, we can prove the commutation relations of a $\phi$-eigenbasis:
\begin{equation}
\left[ e,\phi e\right] =-be+c\phi e+2\xi  \tag{$\pi $}  \label{a}
\end{equation}
\begin{equation}
\left[ \xi ,e\right] =\left( a+\lambda +1\right) \phi e  \tag{$\delta $}
\label{b}
\end{equation}
\begin{equation}
\left[ \phi e,\xi \right] =\left( a-\lambda +1\right) e  \tag{$\tau $}
\label{c}
\end{equation}

\section{An existence theorem\label{example}}

A convenient way to establish the contact system (Definition \ref{cs1}) is through the proof of Theorem \ref{lemma1}. So in this section we discuss this theorem and the contact system will be presented in the next section.

The theorem of Martinet \cite{MA71}, improved by Lutz \cite{LU70}, states that a three dimensional orientable, compact and without boundary manifold admits a contact structure. On the other hand the result of Gromov \cite{G86} states that any non-compact, odd-dimensional manifold carrying a hyperplane field with an almost complex structure accepts a contact structure. The theorem of this section is an existence result concerning both the compact and the non compact case. Moreover the conditions of the theorem are local, hence they can be checked easily.

Let $M^{3}$ be a three dimensional smooth manifold. An absolute parallelism (see e.g. \cite{KNO63} p.122) on $M^{3}$ is a global trivialization of the tangent bundle of $M^{3}$. In other words an absolute parallelism is a $3$-tuple of linearly independent nowhere zero vector fields $\left\{ e_{1},e_{2},e_{3}\right\} $, such that at every $p\in M^{3}$ this $3$-tuple constitutes a basis of $T_{p}M^{3}$. If $M^{3}$ accepts an absolute parallelism then it is said to be \emph{parallelizable}. It is well known that a $3$-dimensional manifold is orientable if and only if it is parallelizable. In what follows $\omega ^{k}$ will denote the dual 1-form of $e_{k}$, $k=1,2,3$. Duality between tangent and cotangent bundle implies $\omega ^{k}\left( e_{j}\right) =\delta _{j}^{k}$.

\begin{theorem}\label{lemma1}
Let $M^{3}$ be a 3-dimensional, orientable, smooth manifold equipped with a Levi-Civita connection $\nabla $. Suppose that

\begin{enumerate}
\item[(1)] Among the trivializations of the tangent bundle of $M^{3}$, there exists one $\left\{ e_{1},e_{2},e_{3}\right\} $ satisfying the following condition: there exists a fixed permutation $i_{0}$, $j_{0}$, $k_{0}$ of $1$, $2$, $3$, such that $f:=\omega ^{k_{0}}\left( \left[ e_{i_{0}},e_{j_{0}}\right] \right) $ is everywhere non zero,

\item[(2)] The integral curves of $e_{k_{0}}$ are geodesics i.e. $\nabla _{e_{k_{0}}}e_{k_{0}}=0$ (same index $k_{0}$ as in (1)),

\item[(3)] The function $f$ is constant along $e_{k_{0}}$ i.e. $e_{k_{0}}\left( f\right) =0$ (same index $k_{0}$ as in (1)).
\end{enumerate}

Then, $M^{3}$ admits a contact metric structure $\left( g,\eta ,\xi ,\phi \right) $. If moreover 
\begin{equation}
e_{j_{0}}\left( f\right) =e_{i_{0}}\left( f\right) =0\text{ and }%
c_{k_{0}j_{0}}^{j_{0}}=c_{k_{0}i_{0}}^{i_{0}}  \label{eigenbasish}
\end{equation}
then $\left\{ e_{1},e_{2},e_{3}\right\} $ is an eigenbasis of $\frac{1}{2} \mathcal{L}_{\xi }\phi $. In this case $c_{k_{0}j_{0}}^{j_{0}}=c_{k_{0}i_{0}}^{i_{0}}=0$ and $f$ is a constant.

Conversely if $M^{3}$ is a contact metric manifold then there exists a global trivialization $\left\{ e_{1}, e_{2}, e_{3}\right\} $ of the tangent bundle of $M^{3}$ satisfying conditions (1), (2) and (3).
\end{theorem}

\begin{proof}
Let $\left\{ e_{1},e_{2},e_{3}\right\} $ be an absolute parallelism of $M^{3} $. The structure functions $c_{ij}^{k}$ are defined by the relation:
\begin{equation*}
\left[ e_{i},e_{j}\right] =c_{ij}^{k}e_{k}
\end{equation*}%
Without causing any damage to generality we may suppose that $\omega^{3}\left( \left[ e_{1},e_{2}\right] \right) =c_{12}^{3}$ is everywhere non zero. For further simplicity we shall denote the function $c_{12}^{3}$ by $f$. If the Euler class $e\left( e_{3}\right) \in H^{2}\left( M,\mathbb{Z}\right) $ of the vector field $e_{3}$ is not zero, then we have to introduce local charts. So on local charts define an almost contact structure on $M^{3}$ as follows: Set $e=e_{1}$, $\phi e=e_{2}$, $\xi =-e_{3}$, and introduce the metric $g\left( e_{i},e_{j}\right) =\delta _{ij}$. Moreover set
\begin{equation}
\phi =\omega ^{1}\otimes e_{2}-\omega ^{2}\otimes e_{1}  \label{u1}
\end{equation}
where $\omega ^{i}$ is the dual vector of $e_{i}$:%
\begin{equation*}
\omega ^{i}\left( e_{j}\right) =\delta _{j}^{i}
\end{equation*}
Then
\begin{equation*}
\left( \omega ^{1}\otimes e_{2}-\omega ^{2}\otimes e_{1}\right) \left(
e_{1}\right) =\omega ^{1}\left( e_{1}\right) e_{2}-\omega ^{2}\left(
e_{1}\right) e_{1}=e_{2}
\end{equation*}
\begin{equation*}
\left( \omega ^{1}\otimes e_{2}-\omega ^{2}\otimes e_{1}\right) \left(
e_{2}\right) =\omega ^{1}\left( e_{2}\right) e_{2}-\omega ^{2}\left(
e_{2}\right) e_{1}=-e_{1}
\end{equation*}
In other words
\begin{equation*}
e_{2}=\phi e_{1},~~~-e_{1}=\phi e_{2},~~~\phi \xi =0
\end{equation*}
Moreover
\begin{equation*}
\Phi =g\left( \cdot ,\phi \cdot \right) =\omega ^{1}\otimes \omega
^{2}-\omega ^{2}\otimes \omega ^{1}
\end{equation*}
It is well known (\cite{B02} p. 53) that $\left( \phi ,\xi ,\omega^{3}\right) $ is an almost contact structure, hence the structure group of the tangent bundle can be reduced to $U\left( 1\right) \times 1$. The transition maps of the tangent bundle respect the metric $\delta _{ij}$ since $\xi $ is globally defined and the symplectic basis $\left\{ e,\phi
e\right\} $ is transformed from chart to chart by the action of an element of $U\left( 1\right) $, since the plane field defined by $\left\{ e,\phi e\right\} $ is also globally defined. Now, from the well known identity $X\omega \left( Y\right) -Y\omega \left( X\right) -\omega \left( \left[ X,Y\right] \right) =d\omega \left( X,Y\right) $ we get
\begin{equation}
d\omega ^{3}\left( e_{1},e_{2}\right) =-\omega ^{3}\left( \left[ e_{1},e_{2}%
\right] \right) =-f  \label{lasd1}
\end{equation}
\begin{equation}
d\omega ^{3}\left( e_{1},e_{3}\right) =-\omega ^{3}\left( \left[ e_{1},e_{3}%
\right] \right) =-c_{13}^{3}  \label{lasd2}
\end{equation}
\begin{equation}
d\omega ^{3}\left( e_{3},e_{2}\right) =-\omega ^{3}\left( \left[ e_{3},e_{2}%
\right] \right) =-c_{32}^{3}  \label{lasd3}
\end{equation}
Since $\nabla g=0$ we have $g\left( \nabla _{e_{i}}e_{j},e_{k}\right)+g\left( e_{j},\nabla _{e_{i}}e_{k}\right) =0$. Putting $i=j=3$ and taking account the hypothesis $\nabla _{e_{3}}e_{3}=0$ we get:
\begin{equation*}
g\left( e_{3},\nabla _{e_{3}}e_{j}\right) =0
\end{equation*}
On the other hand for $j=k=3$ we get $g\left( \nabla_{e_{i}}e_{3},e_{3}\right) +g\left( e_{3},\nabla _{e_{i}}e_{3}\right) =0$. Since $\nabla _{X}Y-\nabla _{Y}X=\left[ X,Y\right] $ we have
\begin{equation}
g\left( \left[ e_{i},e_{3}\right] ,e_{3}\right) =c_{i3}^{3}=0,\text{ }i=1,2,3
\label{conditionb}
\end{equation}
Now, the only non zero component of $\Phi $ is%
\begin{equation}
\Phi \left( e_{1},e_{2}\right) =\omega ^{1}\left( e_{1}\right) \omega
^{2}\left( e_{2}\right) =1  \label{lasd4}
\end{equation}
From (\ref{lasd1}), (\ref{lasd2}), (\ref{lasd3}) and (\ref{lasd4}) we get $-d\omega ^{3}\left( e_{i},e_{j}\right) =f\Phi \left( e_{i},e_{j}\right) $, and thus
\begin{equation*}
\Phi \left( X,Y\right) =-\frac{1}{f}d\omega ^{3}\left( X,Y\right)
\end{equation*}%
for all vector fields $X$, $Y$. Hence the above defined almost contact metric structure $\left( \phi ,e_{3},\omega ^{3},g\right) $ is not an associated one. Note further that $\Phi \wedge \omega ^{3}=2\omega^{1}\wedge \omega ^{2}\wedge \omega ^{3}=-\frac{1}{f}d\omega ^{3}\wedge\omega ^{3}$ is everywhere non zero.

In order to make the almost contact metric structure $\left( \phi ,e_{3},\omega ^{3},g\right) $ an associated one, we apply the deformation
\begin{equation}
\overline{\omega }^{3}:=\frac{f}{2}\omega ^{3},~\overline{e}_{3}:=\frac{2}{f}%
e_{3},~\overline{\phi }:=\phi ,~\overline{g}:=\frac{f^{2}}{4}g
\label{deform1}
\end{equation}
Moreover we set
\begin{equation*}
\overline{e}_{1}:=\frac{2e_{1}}{f}\text{ and }\overline{e}_{2}:=\frac{2e_{1}%
}{f}
\end{equation*}
so that $\overline{g}\left( \overline{e}_{i},\overline{e}_{j}\right) =\delta
_{ij}$. Observe that after the deformation $\overline{\omega }^{3}\wedge d%
\overline{\omega }^{3}=\frac{f^{2}}{4}\omega ^{3}\wedge d\omega ^{3}=\omega
^{1}\wedge \omega ^{2}\wedge \omega ^{3}$, so the positive sign of the
volume form is restored. Now $d\overline{\omega }^{3}=\frac{1}{2}d\left(
f\omega ^{3}\right) =\frac{1}{2}df\wedge \omega ^{3}+\frac{1}{2}fd\omega
^{3} $. From this relation and the identity $\left( a\wedge b\right) \left(
e_{1},e_{2}\right) =a\left( e_{1}\right) b\left( e_{2}\right) -a\left(
e_{2}\right) b\left( e_{1}\right) $ we have%
\begin{equation*}
d\overline{\omega }^{3}\left( e_{1},e_{2}\right) =-\frac{f^{2}}{4}
\end{equation*}%
\begin{equation*}
d\overline{\omega }^{3}\left( e_{1},e_{3}\right) =e_{1}\left( f\right) =0
\end{equation*}%
\begin{equation*}
d\overline{\omega }^{3}\left( e_{2},e_{3}\right) =e_{2}\left( f\right) =0
\end{equation*}%
Hence%
\begin{equation*}
\overline{\Phi }\left( \overline{e}_{1},\overline{e}_{2}\right) =\overline{g}%
\left( \overline{e}_{1},\overline{\phi }\overline{e}_{2}\right) =\frac{4}{%
f^{2}}\frac{f^{2}}{4}g\left( e_{1},\phi e_{2}\right) =-1=d\overline{\omega }%
^{3}\left( \overline{e}_{1},\overline{e}_{2}\right)
\end{equation*}%
\begin{equation*}
\overline{\Phi }\left( \overline{e}_{1},\overline{e}_{3}\right) =\frac{4}{%
f^{2}}\overline{g}\left( e_{1},\overline{\phi }e_{3}\right) =0
\end{equation*}%
\begin{equation*}
\overline{\Phi }\left( \overline{e}_{2},\overline{e}_{3}\right) =\frac{4}{%
f^{2}}\overline{g}\left( e_{3},\overline{\phi }e_{3}\right) =0
\end{equation*}%
Observe that since $\nabla _{\overline{e}_{3}}\overline{e}_{3}=-\frac{4}{%
f^{3}}e_{3}\left( f\right) e_{3}=0$ it follows that $d\overline{\omega }%
^{3}\left( \overline{e}_{1},\overline{e}_{3}\right) =d\overline{\omega }%
^{3}\left( \overline{e}_{2},\overline{e}_{3}\right) =0$. Since $\overline{g}%
\left( \overline{\phi }X,\overline{\phi }Y\right) =\overline{g}\left(
X,Y\right) -\overline{\omega }^{3}\left( X\right) \overline{\omega }%
^{3}\left( Y\right) $ we conclude that $\left( \overline{\phi },\overline{e}%
_{3},\overline{\omega }^{3},\overline{g}\right) $ is an associated contact
metric structure. Note that:

\begin{enumerate}
\item The vector field $\overline{e}_{3}$ is divergence free since $\mathcal{%
L}_{\frac{2}{f}e_{3}}f=0$.

\item The commutation relations of the deformed $\phi $-basis become%
\begin{equation*}
\left[ \overline{e}_{1},\overline{e}_{2}\right] =\frac{4c_{12}^{1}}{f^{2}}%
e_{1}+\frac{4c_{12}^{2}}{f^{2}}e_{2}+\frac{4}{f}e_{3}
\end{equation*}%
\begin{equation*}
\left[ \overline{e}_{2},\overline{e}_{3}\right] =\frac{4c_{23}^{1}}{f^{2}}%
e_{1}+\frac{4c_{23}^{2}}{f^{2}}e_{2}
\end{equation*}%
\begin{equation*}
\left[ \overline{e}_{3},\overline{e}_{1}\right] =\frac{4c_{31}^{1}}{f^{2}}%
e_{1}+\frac{4c_{31}^{2}}{f^{2}}e_{2}
\end{equation*}
\end{enumerate}

A calculation of the Lie derivative of $\phi $ with respect to $\frac{2}{f}%
e_{3}$ yields%
\begin{eqnarray*}
\frac{f}{2}\left( \mathcal{L}_{\frac{2}{f}e_{3}}\phi \right) \left(
e_{i}\right) &=&-c_{3i}^{2}e_{1}+c_{3i}^{1}e_{2}+\left( \delta
_{i}^{2}c_{31}^{1}-\delta _{i}^{1}c_{32}^{1}\right) e_{1}+\left( \delta
_{i}^{2}c_{31}^{2}-\delta _{i}^{1}c_{32}^{2}\right) e_{2}- \\
&&\frac{f}{2}\delta _{i}^{2}e_{1}\left( \frac{2}{f}\right) e_{3}+\frac{f}{2}%
\delta _{i}^{1}e_{2}\left( \frac{2}{f}\right) e_{3}
\end{eqnarray*}%
By definition $2\overline{h}=\mathcal{L}_{\frac{2}{f}e_{3}}\overline{\phi }$%
. Thus the above expression yields%
\begin{equation}
f\overline{h}\left( e_{1}\right) =\left( c_{23}^{1}-c_{31}^{2}\right)
e_{1}+\left( c_{31}^{1}-c_{32}^{2}\right) e_{2}+\frac{f}{2}e_{2}\left( \frac{%
2}{f}\right) e_{3}  \label{eigen1}
\end{equation}%
\begin{equation}
f\overline{h}\left( e_{2}\right) =\left( -c_{32}^{2}+c_{31}^{1}\right)
e_{1}+\left( c_{31}^{2}+c_{32}^{1}\right) e_{2}-\frac{f}{2}e_{1}\left( \frac{%
2}{f}\right) e_{3}  \label{eigen2}
\end{equation}%
\begin{equation}
f\overline{h}\left( e_{3}\right) =0  \label{eigen3}
\end{equation}%
The vector $e_{3}$ is an eigenvector of $\overline{h}$. However the vectors $%
e_{1}$, $e_{2}$ are eigenvectors of $\overline{h}$ if and only if 
\begin{equation}
c_{32}^{2}=c_{31}^{1}\text{ and }e_{1}\left( f\right) =e_{2}\left( f\right)
=0  \label{conditiona}
\end{equation}%
The commutation relations (\ref{eigen1}), (\ref{eigen2}) and (\ref{eigen3})
take the forms%
\begin{equation*}
\left[ \overline{e}_{1},\overline{e}_{2}\right] =\frac{2c_{12}^{1}}{f}%
\overline{e}_{1}+\frac{2c_{12}^{2}}{f}\overline{e}_{2}+2\overline{e}_{3}
\end{equation*}%
\begin{equation*}
\left[ \overline{e}_{2},\overline{e}_{3}\right] =\frac{2c_{23}^{1}}{f}%
\overline{e}_{1}-\frac{2c_{31}^{1}}{f}\overline{e}_{2}
\end{equation*}%
\begin{equation*}
\left[ \overline{e}_{3},\overline{e}_{1}\right] =\frac{2c_{31}^{1}}{f}%
\overline{e}_{1}+\frac{2c_{31}^{2}}{f}\overline{e}_{2}
\end{equation*}%
Now using the fact that on $M^{3}$ $\left( \mathcal{L}_{\phi X}\overline{%
\omega }^{3}\right) \left( Y\right) =\left( \mathcal{L}_{\phi Y}\overline{%
\omega }^{3}\right) \left( X\right) $, $\overline{\Phi }=d\overline{\omega }%
^{3}$, $\nabla _{\overline{e}_{3}}\overline{e}_{3}=0$ and $\overline{g}%
\left( \overline{h}X,Y\right) =\overline{g}\left( X,\overline{h}Y\right) $
(see \cite{B02} p.53) we obtain $\nabla _{X}\overline{e}_{3}=-\overline{\phi 
}X-\overline{\phi }\overline{h}X$. Setting $X=\overline{e}_{1}$ and taking
account that $g\left( \nabla _{\overline{e}_{3}}\overline{e}_{1},\overline{e}%
_{1}\right) =0$ we obtain 
\begin{equation}
c_{31}^{1}=0  \label{conditionc}
\end{equation}%
Thus, the above commutator relations become%
\begin{equation}
\left[ \overline{e}_{1},\overline{e}_{2}\right] =\frac{2c_{12}^{1}}{f}%
\overline{e}_{1}+\frac{2c_{12}^{2}}{f}\overline{e}_{2}+2\overline{e}_{3}
\label{lasteigen1}
\end{equation}%
\begin{equation}
\left[ \overline{e}_{2},\overline{e}_{3}\right] =\frac{2c_{23}^{1}}{f}%
\overline{e}_{1}  \label{lasteigen2}
\end{equation}%
\begin{equation}
\left[ \overline{e}_{3},\overline{e}_{1}\right] =\frac{2c_{31}^{2}}{f}%
\overline{e}_{2}  \label{lasteigen3}
\end{equation}%
Setting $\frac{2c_{23}^{1}}{f}=a-\lambda +1$, $\frac{2c_{31}^{2}}{f}%
=a+\lambda +1$,%
\begin{equation}
\frac{2c_{12}^{1}}{f}=-b,\text{ and }\frac{2c_{12}^{2}}{f}=c  \label{bc}
\end{equation}%
we obtain the relations ($\pi $), ($\delta $) and ($\tau $) with 
\begin{equation}
\lambda =\frac{c_{31}^{2}-c_{23}^{1}}{f}\text{ and }a=-1+\frac{%
c_{31}^{2}+c_{23}^{1}}{f}  \label{al}
\end{equation}

Conversely let $(M^{3},\phi ,g,\eta ,\xi )$ be a contact metric manifold.
Since $M^{3}$ is orientable it is parallelizable. Choose a parallelization $%
\left\{ e_{1},e_{2},e_{3}\right\} $ with $e_{3}=\xi $. Consider the sets%
\begin{equation*}
U_{0}=\left\{ m\in M^{3}:\lambda =0\text{ in a neighborhood of }m\right\}
\end{equation*}%
\begin{equation*}
U_{1}=\left\{ m\in M^{3}:\lambda \neq 0\text{ in a neighborhood of }m\right\}
\end{equation*}%
On $U_{1}$ we can take $e_{1}$ and $e_{2}$ as eigenvectors of $h$ with
eigenvalues $\lambda $ and $-\lambda $ respectively. Then the conditions 
\emph{(1)}, \emph{(2)} and \emph{(3)} of the present Theorem are fulfilled
since by hypothesis $M^{3}$ is a contact metric manifold hence $f=2$ and $%
\nabla _{e_{3}}e_{3}=0$.

$U_{0}$ is Sasakian. Hence $\nabla _{e_{3}}e_{3}=0$ and thus (2) is
fulfilled. Condition (1) also holds due to the non-integrability of the
contact subbundle $D$ (if $\left[ e_{1},e_{2}\right] \in D$ then $D$ would
be integrable). Condition (3) is also true: From identity%
\begin{equation}
\left( \mathcal{L}_{\xi X}g\right) \left( Y_{1},Y_{2}\right) =X\cdot g\left(
Y_{1},Y_{2}\right) -g\left( \left[ X,Y_{1}\right] ,Y_{2}\right) -g\left(
Y_{1},\left[ X,Y_{2}\right] \right)  \label{id1}
\end{equation}%
setting $X=\xi ,Y_{1}=\left[ e,\phi e\right] ,Y_{2}=\xi $ we get the relation $\left( 
\mathcal{L}_{\xi }g\right) \left( \left[ e,\phi e\right] ,\xi \right) =\xi
\cdot g\left( \left[ e,\phi e\right] ,\xi \right) -g\left( \left[ \xi ,\left[
e,\phi e\right] \right] ,\xi \right) -g\left( \left[ e,\phi e\right] ,\left[
\xi ,\xi \right] \right) $. Since $\mathcal{L}_{\xi }g=0$ and $g\left( \left[
e,\phi e\right] ,\left[ \xi ,\xi \right] \right) =0$ the identity (\ref{id1}%
) gives $\xi \cdot g\left( \left[ e,\phi e\right] ,\xi \right) =g\left( %
\left[ \xi ,\left[ e,\phi e\right] \right] ,\xi \right) $. Using Lemma \ref%
{aux1} and after some calculations we get the relation $g\left( \left[ \xi ,\left[ e,\phi
e\right] \right] ,\xi \right) =cg\left( \left[ \xi ,\phi e\right] ,\xi
\right) -bg\left( \left[ \xi ,e\right] ,\xi \right)$, which is zero. Hence $%
\xi \cdot g\left( \left[ e,\phi e\right] ,\xi \right) =0$.

Properties (1), (2) and (3) are true on $U_{0}\cup U_{1}$ which is open and
dense in $M^{3}$. But these properties involve equations between continuous
functions. Hence they are true all over $M^{3}$.
\end{proof}

\begin{remark}
The following question arises: Is the definition of $\phi $ by (\ref{u1})
and the transformation (\ref{deform1}) unique? In other words, are there any
other definitions for $\phi $ and transformations of the resulting almost
contact structure which give rise to a contact metric structure on $M^{3}$?
If no then the $\phi $ defined by (\ref{u1}) and the transformation (\ref%
{deform1}) are universal. If yes then is there any connection among the
resulting contact structures? Note that in the later case, if the resulting
contact structures constitute a $1$-parameter smooth family and if the
manifold is closed then by a Theorem of Gray \cite{GR59} they are all
isotopic.
\end{remark}

\section{The contact system}\label{par4456}

\subsection{The form of the contact system}

Let $M^{3}$ be an orientable $3$-manifold. Choose on $M^{3}$ a local chart $%
\left( U^{\prime },x\right) $ and on $U:=x\left( U^{\prime }\right) \subset 
\mathbb{R}^{3}$ introduce a positively oriented coordinate basis $\left\{ 
\frac{\partial }{\partial x^{1}},\frac{\partial }{\partial x^{2}},\frac{\partial }{\partial x^{3}}\right\}$. If $g$ is a Riemannian metric on $%
M^{3} $ we have $g_{ij}=g\left( \frac{\partial }{\partial x^{i}},\frac{%
\partial }{\partial x^{j}}\right) $. Let $\left\{ e_{1},e_{2},e_{3}\right\} $
be a positively oriented orthonormal parallelization of $U$, i.e. $\delta
_{ab}=g\left( e_{\alpha },e_{\beta }\right) $. On $U$ there exist nine
functions $b_{\alpha }^{i}$, $i,\alpha =1,2,3$ such that%
\begin{equation}
\left[ 
\begin{array}{c}
e_{1} \\ 
e_{2} \\ 
e_{3}%
\end{array}%
\right] =\left[ 
\begin{array}{ccc}
b_{1}^{1} & b_{1}^{2} & b_{1}^{3} \\ 
b_{2}^{1} & b_{2}^{2} & b_{2}^{3} \\ 
b_{3}^{1} & b_{3}^{2} & b_{3}^{3}%
\end{array}%
\right] \left[ 
\begin{array}{c}
\frac{\partial }{\partial x^{1}} \\ 
\frac{\partial }{\partial x^{2}} \\ 
\frac{\partial }{\partial x^{3}}%
\end{array}%
\right] ,\text{ and }\left[ 
\begin{array}{c}
\frac{\partial }{\partial x^{1}} \\ 
\frac{\partial }{\partial x^{2}} \\ 
\frac{\partial }{\partial x^{3}}%
\end{array}%
\right] =\left[ 
\begin{array}{ccc}
b_{1}^{1} & b_{1}^{2} & b_{1}^{3} \\ 
b_{2}^{1} & b_{2}^{2} & b_{2}^{3} \\ 
b_{3}^{1} & b_{3}^{2} & b_{3}^{3}%
\end{array}%
\right] ^{-1}\left[ 
\begin{array}{c}
e_{1} \\ 
e_{2} \\ 
e_{3}%
\end{array}%
\right]  \label{matrixform1}
\end{equation}%
where $b_{\alpha }^{i}\left( b^{-1}\right) _{i}^{\beta }=\delta _{\alpha
}^{\beta }$ and $\left( b^{-1}\right) _{i}^{\alpha }b_{\alpha }^{j}=\delta
_{i}^{j}$.

\begin{definition}
We call $B$-matrix the matrix whose elements are the functions $b_{\alpha
}^{i}$ defined in (\ref{matrixform1}). The function $\left( b^{-1}\right)
_{i}^{\alpha }$, $\alpha =1,2,3$, $i=1,2,3$, is the element of the $i$ row
and the $\alpha $ column of the matrix $B^{-1}$. It is obvious that $\det
B>0 $ everywhere.
\end{definition}

We seek a necessary and sufficient condition such that the basis $\left\{
e_{1}, e_{2}, e_{3}\right\}$ in (\ref{matrixform1}) to be a $\phi$-eigenbasis of a contact metric structure $\left( g,\eta ,\xi ,\phi \right)$ of $M^{3}$. In terms of the structure functions $c_{ij}^{k}$ these conditions are already found in the course of the proof of Theorem \ref{lemma1}. These are (\ref{conditionb}), (\ref{conditiona}) and (\ref{conditionc}). Let us summarize these conditions:%
\begin{equation}
c_{12}^{3}=2,\text{ }c_{31}^{3}=c_{32}^{3}=0,\text{ }c_{32}^{2}=c_{31}^{1}=0
\label{easy1}
\end{equation}%
Define the functions $c_{\alpha \beta }^{\gamma }$ by the relation $\left[
e_{\alpha },e_{\beta }\right] =c_{\alpha \beta }^{\gamma }e_{\gamma }$. Then
using (\ref{matrixform1}) we obtain%
\begin{equation}
c_{\alpha \beta }^{\gamma }=\left( b^{-1}\right) _{j}^{\gamma }\left(
b_{\alpha }^{i}\frac{\partial b_{\beta }^{j}}{\partial x^{i}}-b_{\beta }^{i}%
\frac{\partial b_{\alpha }^{j}}{\partial x^{i}}\right)  \label{struct1}
\end{equation}%
where $x^{i}$, $i=1,2,3$ are local coordinates of $U$. Calculating $B^{-1}$
and using (\ref{easy1}) and (\ref{struct1}) we have the following

\begin{theorem}
The following equations are necessary and sufficient conditions in order for
the basis $\left\{ e_{1},e_{2},e_{3}\right\} $ to be a $\phi $-eigenbasis of
a contact metric structure $\left( g,\eta ,\xi ,\phi \right) $ of $M^{3}$.%
\begin{gather}
\left( b_{1}^{2}b_{2}^{3}-b_{1}^{3}b_{2}^{2}\right) \left( b_{3}^{i}\frac{%
\partial b_{1}^{1}}{\partial x^{i}}-b_{1}^{i}\frac{\partial b_{3}^{1}}{%
\partial x^{i}}\right) +  \label{contact1} \\
\left( b_{1}^{3}b_{2}^{1}-b_{1}^{1}b_{2}^{3}\right) \left( b_{3}^{i}\frac{%
\partial b_{1}^{2}}{\partial x^{i}}-b_{1}^{i}\frac{\partial b_{3}^{2}}{%
\partial x^{i}}\right) +  \notag \\
\left( b_{1}^{1}b_{2}^{2}-b_{1}^{2}b_{2}^{1}\right) \left( b_{3}^{i}\frac{%
\partial b_{1}^{3}}{\partial x^{i}}-b_{1}^{i}\frac{\partial b_{3}^{3}}{%
\partial x^{i}}\right) =0  \notag
\end{gather}%
\begin{gather}
\left( b_{1}^{2}b_{2}^{3}-b_{1}^{3}b_{2}^{2}\right) \left( b_{3}^{i}\frac{%
\partial b_{2}^{1}}{\partial x^{i}}-b_{2}^{i}\frac{\partial b_{3}^{1}}{%
\partial x^{i}}\right) +  \label{contact2} \\
\left( b_{1}^{3}b_{2}^{1}-b_{1}^{1}b_{2}^{3}\right) \left( b_{3}^{i}\frac{%
\partial b_{2}^{2}}{\partial x^{i}}-b_{2}^{i}\frac{\partial b_{3}^{2}}{%
\partial x^{i}}\right) +  \notag \\
\left( b_{1}^{1}b_{2}^{2}-b_{1}^{2}b_{2}^{1}\right) \left( b_{3}^{i}\frac{%
\partial b_{2}^{3}}{\partial x^{i}}-b_{2}^{i}\frac{\partial b_{3}^{3}}{%
\partial x^{i}}\right) =0  \notag
\end{gather}%
\begin{gather}
\left( b_{1}^{2}b_{2}^{3}-b_{1}^{3}b_{2}^{2}\right) \left( b_{1}^{i}\frac{%
\partial b_{2}^{1}}{\partial x^{i}}-b_{2}^{i}\frac{\partial b_{1}^{1}}{%
\partial x^{i}}\right) +  \label{contact3} \\
\left( b_{1}^{3}b_{2}^{1}-b_{1}^{1}b_{2}^{3}\right) \left( b_{1}^{i}\frac{%
\partial b_{2}^{2}}{\partial x^{i}}-b_{2}^{i}\frac{\partial b_{1}^{2}}{%
\partial x^{i}}\right) +  \notag \\
\left( b_{1}^{1}b_{2}^{2}-b_{1}^{2}b_{2}^{1}\right) \left( b_{1}^{i}\frac{%
\partial b_{2}^{3}}{\partial x^{i}}-b_{2}^{i}\frac{\partial b_{1}^{3}}{%
\partial x^{i}}\right) =2\det B  \notag
\end{gather}%
\begin{gather}
\left( b_{2}^{2}b_{3}^{3}-b_{2}^{3}b_{3}^{2}\right) \left( b_{3}^{i}\frac{%
\partial b_{1}^{1}}{\partial x^{i}}-b_{1}^{i}\frac{\partial b_{3}^{1}}{%
\partial x^{i}}\right) +  \label{contact4} \\
\left( b_{2}^{3}b_{3}^{1}-b_{2}^{1}b_{3}^{3}\right) \left( b_{3}^{i}\frac{%
\partial b_{1}^{2}}{\partial x^{i}}-b_{1}^{i}\frac{\partial b_{3}^{2}}{%
\partial x^{i}}\right) +  \notag \\
\left( b_{2}^{1}b_{3}^{2}-b_{2}^{2}b_{3}^{1}\right) \left( b_{3}^{i}\frac{%
\partial b_{1}^{3}}{\partial x^{i}}-b_{1}^{i}\frac{\partial b_{3}^{3}}{%
\partial x^{i}}\right) =0  \notag
\end{gather}%
\begin{gather}
\left( b_{1}^{3}b_{3}^{2}-b_{1}^{2}b_{3}^{3}\right) \left( b_{3}^{i}\frac{%
\partial b_{2}^{1}}{\partial x^{i}}-b_{2}^{i}\frac{\partial b_{3}^{1}}{%
\partial x^{i}}\right) +  \label{contact5} \\
\left( b_{1}^{1}b_{3}^{3}-b_{1}^{3}b_{3}^{1}\right) \left( b_{3}^{i}\frac{%
\partial b_{2}^{2}}{\partial x^{i}}-b_{2}^{i}\frac{\partial b_{3}^{2}}{%
\partial x^{i}}\right) +  \notag \\
\left( b_{1}^{2}b_{3}^{1}-b_{1}^{1}b_{3}^{2}\right) \left( b_{3}^{i}\frac{%
\partial b_{2}^{3}}{\partial x^{i}}-b_{2}^{i}\frac{\partial b_{3}^{3}}{%
\partial x^{i}}\right) =0  \notag
\end{gather}
\end{theorem}

\begin{definition}
\label{cs1}Equations (\ref{contact1}) - (\ref{contact5}) are called the
contact system.
\end{definition}

\begin{remark}
The metric of $M^{3}$ in local coordinates $x$ is given by the relation%
\begin{equation}
g_{ij}=\left( b^{-1}\right) _{i}^{m}\left( b^{-1}\right) _{j}^{n}\delta _{mn}
\label{m2}
\end{equation}
\end{remark}

\section{Reduction and solution of the contact system\label{sectsim}}

\subsection{A Lemma}

By the following Lemma, we will simplify the $B$-matrix:

\begin{lemma}
\label{lemmaauxsimple}Let $M^{3}$ be an orientable $3$-manifold and let $%
\left( U^{\prime },\theta \right) $ be a local chart of $M^{3}$. Let $%
\left\{ e_{1},e_{2},e_{3}\right\} $ be an orthonormal parallelization of $%
M^{3}$ and let $U=\theta \left( U^{\prime }\right) \subset \mathbb{R}^{3}$.
There exist on $U$ local coordinates $\omega $ such that%
\begin{equation*}
e_{3}=\frac{\partial }{\partial \omega ^{1}}
\end{equation*}%
and%
\begin{equation*}
e_{1}=\chi _{1}\frac{\partial }{\partial \omega ^{1}}+\chi _{2}\frac{%
\partial }{\partial \omega ^{2}}
\end{equation*}%
where $\chi _{1}$ and $\chi _{2}$ are smooth functions on $U$.
\end{lemma}

\begin{proof}
We may write on $U$:%
\begin{equation}
e_{i}=b_{i}^{\alpha }\left( \theta \right) \frac{\partial }{\partial \theta
^{\alpha }}=b_{i}^{\alpha }\left( \theta \right) \frac{\partial \psi ^{\beta
}}{\partial \theta ^{\alpha }}\frac{\partial }{\partial \psi ^{\beta }}
\label{gold1}
\end{equation}%
where $\psi $ are some other coordinates on $U$. We claim that $\psi $ can
be chosen so that the following relations hold:%
\begin{equation}
b_{3}^{\alpha }\left( \theta \right) \frac{\partial \psi ^{1}}{\partial
\theta ^{\alpha }}=1,\text{ }b_{3}^{\alpha }\left( \theta \right) \frac{%
\partial \psi ^{2}}{\partial \theta ^{\alpha }}=0,\text{ }b_{3}^{\alpha
}\left( \theta \right) \frac{\partial \psi ^{3}}{\partial \theta ^{\alpha }}%
=0,~b_{1}^{\alpha }\left( \theta \right) \frac{\partial \psi ^{3}}{\partial
\theta ^{\alpha }}=0  \label{lemmaaux1}
\end{equation}%
Set%
\begin{equation}
\frac{\partial \psi ^{1}}{\partial \theta ^{\beta }}:=g_{\beta \gamma
}\left( \theta \right) b_{3}^{\gamma }\left( \theta \right) ,\text{ }\frac{%
\partial \psi ^{2}}{\partial \theta ^{\beta }}:=g_{\beta \gamma }\left(
\theta \right) b_{1}^{\gamma }\left( \theta \right) ,\text{ }\frac{\partial
\psi ^{3}}{\partial \theta ^{\beta }}:=g_{\beta \gamma }\left( \theta
\right) b_{2}^{\gamma }\left( \theta \right)  \label{lemmaaux2}
\end{equation}%
and observe that, since $g\left( e_{3},e_{3}\right) =1$, we have%
\begin{equation}
g_{\alpha \gamma }\left( \theta \right) b_{3}^{\gamma }\left( \theta \right)
b_{3}^{\alpha }\left( \theta \right) =1  \label{lemmaaux3}
\end{equation}%
Since $g\left( e_{2},e_{3}\right) =g\left( e_{1},e_{3}\right) =g\left(
e_{1},e_{2}\right) =0$ we must have%
\begin{equation}
g_{\alpha \gamma }\left( \theta \right) b_{3}^{\gamma }\left( \theta \right)
b_{1}^{\alpha }\left( \theta \right) =g_{\alpha \gamma }\left( \theta
\right) b_{3}^{\gamma }\left( \theta \right) b_{2}^{\alpha }\left( \theta
\right) =g_{\alpha \gamma }\left( \theta \right) b_{1}^{\gamma }\left(
\theta \right) b_{2}^{\alpha }\left( \theta \right) =0  \label{lemmaaux4}
\end{equation}%
In view of (\ref{lemmaaux3}) and (\ref{lemmaaux4}) we see that the
coordinates $\left( \psi ^{1},\psi ^{2},\psi ^{3}\right) $ defined by (\ref%
{lemmaaux2}) satisfy the claimed relations, namely (\ref{lemmaaux1}). Then
using (\ref{gold1}) and setting $\chi _{1}=b_{1}^{\alpha }\left( \theta
\right) \frac{\partial \psi ^{1}}{\partial \theta ^{\alpha }}$, $\chi
_{2}=b_{1}^{\alpha }\left( \theta \right) \frac{\partial \psi ^{2}}{\partial
\theta ^{\alpha }}$ and $\omega =\psi $, the Lemma is proved.
\end{proof}

\subsection{The reduction}

From Lemma \ref{lemmaauxsimple} we see that there exists a coordinate system
on $U$ such that the $B$-matrix takes the form%
\begin{equation}
B=\left[ 
\begin{array}{ccc}
b_{1}^{1} & b_{1}^{2} & 0 \\ 
b_{2}^{1} & b_{2}^{2} & b_{2}^{3} \\ 
1 & 0 & 0%
\end{array}%
\right]  \label{simplified1}
\end{equation}

\begin{definition}
\label{simpdef}A local coordinate system on which the $B$-matrix takes the
form (\ref{simplified1}) is called simplifying coordinates. From now on we
will set%
\begin{equation}
\left[ 
\begin{array}{ccc}
b_{1}^{1} & b_{1}^{2} & 0 \\ 
b_{2}^{1} & b_{2}^{2} & b_{2}^{3} \\ 
1 & 0 & 0%
\end{array}%
\right] :=\left[ 
\begin{array}{ccc}
\alpha & \beta & 0 \\ 
\delta & \epsilon & \zeta \\ 
1 & 0 & 0%
\end{array}%
\right]  \label{simplified2}
\end{equation}
\end{definition}

Taking account (\ref{simplified1}) as well as (\ref{contact1}) - (\ref%
{contact5}), and the fact that $\det B=\beta \zeta \neq 0$, the contact
system becomes in simplifying coordinates:%
\begin{equation}
\alpha _{1}=0  \label{S1}
\end{equation}%
\begin{equation}
\beta \delta _{1}=\alpha \epsilon _{1}  \label{S2}
\end{equation}%
\begin{equation}
\beta \left( \alpha \epsilon -\delta \beta \right) \zeta _{2}+\left( \alpha
\beta _{3}-\beta \alpha _{3}\right) \zeta ^{2}+\left( \beta ^{2}\delta
_{2}-\alpha \beta \epsilon _{2}+\alpha \epsilon \beta _{2}-\beta \epsilon
\alpha _{2}-2\beta \right) \zeta =0  \label{S3}
\end{equation}%
\begin{equation}
\beta _{1}=0  \label{S4}
\end{equation}%
\begin{equation}
\zeta _{1}=0  \label{S5}
\end{equation}

The relation (\ref{S2}) can be written as $\delta _{1}=\frac{\alpha \epsilon
_{1}}{\beta }$. Using (\ref{S1}) and (\ref{S4}) we get%
\begin{equation}
\delta =\frac{\alpha \epsilon }{\beta }+F\left( x^{2},x^{3}\right)
\label{T2}
\end{equation}%
where $F\left( x^{2},x^{3}\right) $ is an arbitrary smooth function. Then $\delta _{2}$ equals to $F_{2}+\frac{\alpha _{2}\epsilon \beta +\alpha \epsilon _{2}\beta
-\alpha \epsilon \beta _{2}}{\beta ^{2}}$ and (\ref{S3}) becomes%
\begin{equation}
F\zeta _{2}+\left( \frac{\alpha }{\beta }\right) _{3}\zeta ^{2}-\left( F_{2}-%
\frac{2}{\beta }\right) \zeta =0  \label{T3}
\end{equation}%
If $F=0$ then (\ref{T3}) becomes a linear equation with respect to $\zeta $,
namely%
\begin{equation}
\left( \frac{\alpha }{\beta }\right) _{3}\zeta =\left( F_{2}-\frac{2}{\beta }%
\right)  \label{T*3}
\end{equation}%
If $F_{2}=0$ and $\left( \frac{\alpha }{\beta }\right) _{3}=0$ then $\beta
=0 $ which is absurd, since by hypothesis $\beta \neq 0$. If $F_{2}\neq 0$
and $\left( \frac{\alpha }{\beta }\right) _{3}=0$ then $F_{2}=\frac{2}{\beta 
}$ and $\zeta $ is any non-zero function of $x^{2}$ and $x^{3}$. If $\left( 
\frac{\alpha }{\beta }\right) _{3}\neq 0$ then $\zeta =-\frac{2}{\beta
\left( \frac{\alpha }{\beta }\right) _{3}}$. If $F\neq 0$ then (\ref{T3})
yields%
\begin{equation}
\zeta _{2}+\frac{1}{F}\left( \frac{\alpha }{\beta }\right) _{3}\zeta
^{2}-\left( \frac{F_{2}}{F}-\frac{2}{F\beta }\right) \zeta =0  \label{S*3}
\end{equation}%
The equation (\ref{S*3}) is a Riccati ordinary differential equation. The
unknown function is $\zeta \left( x^{2},x^{3}\right) $ and the variable is $%
x^{2}$ while $x^{3}$ is treated as a parameter. The solution of the above
equation is%
\begin{equation}
\zeta =F\frac{e^{-2\int \frac{dx^{2}}{\beta F}}}{\int e^{-2\int \frac{dx^{2}%
}{\beta F}}\left( \frac{\alpha }{\beta }\right) _{3}dx^{2}-K\left(
x^{3}\right) }  \label{sur}
\end{equation}%
where $K\left( x^{3}\right) $ is a smooth function. Observe that the right
hand side of (\ref{sur}) is a function of $x^{2}$ and $x^{3}$ in accordance
with (\ref{S5}).

\subsection{The contact form and the metric}

Define the $1$-form $\eta $ by $\eta \left( X\right) =g\left( X,e_{3}\right) 
$. In simplifying coordinates we have $e_{3}=\frac{\partial }{\partial x^{3}}
$ and $\eta =dx^{1}-\frac{\alpha }{\beta }dx^{2}-\frac{F}{\zeta }dx^{3}$. It
is easy to see that (\ref{S3}) can be written as%
\begin{equation*}
\left[ \frac{\alpha \epsilon -\delta \beta }{\beta \zeta }\right]
_{2}+\left( \frac{\alpha }{\beta }\right) _{3}=-\frac{2}{\beta \zeta }
\end{equation*}%
A calculation shows that%
\begin{equation*}
\eta \wedge d\eta =\frac{1}{2}\left[ -\left( \frac{\alpha }{\beta }\right)
_{3}-\left( \frac{\alpha \epsilon -\beta \delta }{\beta \zeta }\right) _{2}%
\right] dx^{1}\wedge dx^{2}\wedge dx^{3}
\end{equation*}%
Hence $\eta \wedge d\eta =\frac{1}{\beta \zeta }dx^{1}\wedge dx^{2}\wedge
dx^{3}\neq 0$ everywhere, and thus $\eta $ is a contact form, as expected.

The formula (\ref{m2}) yields%
\begin{equation}
g=\left[ 
\begin{array}{ccc}
1 & -\frac{\alpha }{\beta } & \frac{-F}{\zeta } \\ 
-\frac{\alpha }{\beta } & \frac{1+\alpha ^{2}}{\beta ^{2}} & \frac{\alpha
\beta F-\epsilon }{\beta ^{2}\zeta } \\ 
\frac{-F}{\zeta } & \frac{\alpha \beta F-\epsilon }{\beta ^{2}\zeta } & 
\frac{\beta ^{2}\left( 1+F^{2}\right) +\epsilon ^{2}}{\beta ^{2}\zeta ^{2}}%
\end{array}%
\right]  \label{metric-simp}
\end{equation}%
and $\det g=\frac{1}{\beta ^{2}\zeta ^{2}}$. The relation (\ref{metric-simp}%
) along with (\ref{sur}), (\ref{T*3}), (\ref{S4}) and (\ref{S5}) yield the
local classification of $3$-dimensional contact manifolds.

\section{Construction of contact metric structures on $3$-manifolds}

\subsection{An algorithm for the construction of contact metric structures
on $3$-manifolds\label{algorithm}}

Summarizing the results of the previous section, we give an algorithm
producing the contact metric structures on $3$-manifolds.

\begin{enumerate}
\item Consider a $3$-dimensional orientable manifold $M^{3}$ with a
parallelization $\left\{ e_{1},e_{2},e_{3}\right\} $, a local
chart $\left( U^{\prime },\theta \right) $ and the metric $g\left(
e_{i},e_{j}\right) =\delta _{ij}$. Consider the open $\theta \left( U\right)
\subset \mathbb{R}^{3}$. Introduce on $\theta \left( U\right) $ simplifying
coordinates, $x^{1}$, $x^{2}$ and $x^{3}$.

\item \label{step2}Insert functions $\alpha \left( x^{2},x^{3}\right) $, $%
\beta \left( x^{2},x^{3}\right) \neq 0$, $\epsilon \left(
x^{1},x^{2},x^{3}\right) $ and $F\left( x^{2},x^{3}\right) $.\\ Set $\delta
\left( x^{1},x^{2},x^{3}\right) :=\frac{\alpha \epsilon }{\beta }+F$.

\item If $F=0$ then compute the quantity $\left( \frac{\alpha }{\beta }%
\right) _{3}$.

\begin{enumerate}
\item If $\left( \frac{\alpha }{\beta }\right) _{3}=0$ then examine if $%
F_{2}=\frac{2}{\beta }$.

\begin{enumerate}
\item If $F_{2}=\frac{2}{\beta }$ holds then set $\zeta $ to be any non-zero
smooth function of $x^{2}$ and $x^{3}$ and go to step \ref{step}.

\item If $F_{2}\neq \frac{2}{\beta }$ then go to step \ref{step2} and try an 
$F$ for which $F_{2}=\frac{2}{\beta }$.
\end{enumerate}

\item If $\left( \frac{\alpha }{\beta }\right) _{3}\neq 0$ then set $\zeta =%
\frac{F_{2}-\frac{2}{\beta }}{\left( \frac{\alpha }{\beta }\right) _{3}}$
and go to step \ref{step}.
\end{enumerate}

\item If $F\neq 0$ then compute $\zeta $ by the relation%
\begin{equation*}
\zeta =F\frac{e^{-2\int \frac{dx^{2}}{\beta F}}}{\int e^{-2\int \frac{dx^{2}%
}{\beta F}}\left( \frac{\alpha }{\beta }\right) _{3}dx^{2}-K\left(
x^{3}\right) }
\end{equation*}

\item \label{step}Compute the metric $g$ in simplifying coordinates by the
relation:%
\begin{equation*}
g=\left[ 
\begin{array}{ccc}
1 & -\frac{\alpha }{\beta } & \frac{-F}{\zeta } \\ 
-\frac{\alpha }{\beta } & \frac{1+\alpha ^{2}}{\beta ^{2}} & \frac{\alpha
\beta F-\epsilon }{\beta ^{2}\zeta } \\ 
\frac{-F}{\zeta } & \frac{\alpha \beta F-\epsilon }{\beta ^{2}\zeta } & 
\frac{\beta ^{2}\left( 1+F^{2}\right) +\epsilon ^{2}}{\beta ^{2}\zeta ^{2}}%
\end{array}%
\right]
\end{equation*}

\item Write down the $B$-matrix%
\begin{equation*}
B=\left[ 
\begin{array}{ccc}
\alpha & \beta & 0 \\ 
\delta & \epsilon & \zeta \\ 
1 & 0 & 0%
\end{array}%
\right]
\end{equation*}
Using Lemma \ref{lemmaauxsimple} write down the vector fields $e_{1}=\alpha 
\frac{\partial }{\partial x^{1}}+\beta \frac{\partial }{\partial x^{2}}$, $%
e_{2}=\delta \frac{\partial }{\partial x^{1}}+\epsilon \frac{\partial }{%
\partial x^{2}}+\zeta \frac{\partial }{\partial x^{3}}$ and $e_{3}=\frac{%
\partial }{\partial x^{3}}$ on $U$ in simplifying coordinates.

\item Define on $M^{3}$ the the $1$-form $\eta $ by $\eta \left( X\right)
=g\left( X,e_{3}\right) $. On $U$ and in simplifying coordinates we have $%
\eta =dx^{1}-\frac{\alpha }{\beta }dx^{2}-\frac{F}{\zeta }dx^{3}$. A
calculation shows that $\eta \wedge d\eta =\frac{1}{\beta \zeta }%
dx^{1}\wedge dx^{2}\wedge dx^{3}\neq 0$ everywhere on $U$, thus $\eta $ is a
contact form.

\item Define the $(1,1)$ tensor field $\phi $ by $\phi e_{3}=0$, $\phi
e_{1}=e_{2}$, $\phi e_{2}=-e_{1}$. Then $\eta \left( e_{3}\right) =1$, $\phi
^{2}X=-X+\eta \left( X\right) e_{3}$, $d\eta \left( X,Y\right) =g\left(
X,\phi Y\right) $ and $g\left( \phi X,\phi Y\right) =g\left( X,Y\right)
-\eta \left( X\right) \eta \left( Y\right) $ for all vector fields $X$, $Y$
on $M^{3}$.

\item Put $\xi :=e_{3}$, $\phi e:=e_{2}$, $e:=e_{1}$. The tetrad $(\eta $, $%
\phi $, $\xi $, $g)$ is a contact metric structure.
\end{enumerate}

\subsection{Example\label{example1}}

We give an example of a contact metric $3$-manifold. We use the algorithm of
section \ref{algorithm} and the notion of simplifying coordinates.

Consider $M^{3}=\mathbb{R}^{3}-\left\{ x_{2}=0\right\} $ with simplifying
coordinates $x^{1}$, $x^{2}$ and $x^{3}$. Let $\alpha \left(
x^{2},x^{3}\right) =0$, $\beta \left( x^{2},x^{3}\right) =x^{2}$, $\epsilon
\left( x^{1},x^{2},x^{3}\right) =x^{1}$ and $F\left( x^{2},x^{3}\right) =1$.
Set $\delta \left( x^{1},x^{2},x^{3}\right) :=\frac{\alpha \epsilon }{\beta }%
+F=1$. Since $F\neq 0$ we have $\zeta =-\frac{1}{\left( x^{2}\right)
^{2}K\left( x^{3}\right) }$. The metric is%
\begin{equation*}
g=\left[ 
\begin{array}{ccc}
1 & 0 & \left( x^{2}\right) ^{2}K\left( x^{3}\right)  \\ 
0 & \frac{1}{\left( x^{2}\right) ^{2}} & x^{1}K\left( x^{3}\right)  \\ 
\left( x^{2}\right) ^{2}K\left( x^{3}\right)  & x^{1}K\left( x^{3}\right)  & 
\left( x^{2}\right) ^{2}K^{2}\left( x^{3}\right) \left( \left( x^{1}\right)
^{2}+2\left( x^{2}\right) ^{2}\right) 
\end{array}%
\right] 
\end{equation*}%
The $B$-matrix is%
\begin{equation*}
B=\left[ 
\begin{array}{ccc}
0 & x^{2} & 0 \\ 
1 & x^{1} & -\frac{1}{\left( x^{2}\right) ^{2}K\left( x^{3}\right) } \\ 
1 & 0 & 0%
\end{array}%
\right] 
\end{equation*}%
and define the vector fields $e_{1}:=x^{2}\frac{\partial }{\partial x^{2}}$, 
$e_{2}:=\frac{\partial }{\partial x^{1}}+x^{1}\frac{\partial }{\partial x^{2}%
}-\frac{1}{\left( x^{2}\right) ^{2}K\left( x^{3}\right) }\frac{\partial }{%
\partial x^{3}}$ and $e_{3}:=\frac{\partial }{\partial x^{3}}$. Define the
the $1$-form $\eta $ by $\eta \left( X\right) =g\left( X,e_{3}\right) $. In
simplifying coordinates we have $\eta =dx^{1}+\left( x^{2}\right)
^{2}K\left( x^{3}\right) dx^{3}$ and $\eta \wedge d\eta =-\frac{1}{%
x^{2}K\left( x^{3}\right) }dx^{1}\wedge dx^{2}\wedge dx^{3}\neq 0$
everywhere. Define the $(1,1)$ tensor field $\phi $ by $\phi e_{3}=0$, $\phi
e_{1}=e_{2}$, $\phi e_{2}=-e_{1}$. Then $\eta \left( e_{3}\right) =1$, $\phi
^{2}X=-X+\eta \left( X\right) e_{3}$, $d\eta \left( X,Y\right) =g\left(
X,\phi Y\right) $ and $g\left( \phi X,\phi Y\right) =g\left( X,Y\right)
-\eta \left( X\right) \eta \left( Y\right) $ for all vector fields $X$, $Y$
on $M^{3}$. Put $\xi :=e_{3}$, $\phi e:=e_{2}$, $e:=e_{1}$. The tetrad $%
(\eta $, $\phi $, $\xi $, $g)$ is a contact metric structure on $M^{3}$.
Further calculations yield%
\begin{equation*}
\left[ e_{1},e_{2}\right] =-2e_{1}-\frac{x^{1}}{x^{2}}e_{2}+2e_{3},~~~\left[
e_{2},e_{3}\right] =-\frac{1}{x^{2}}e_{1},~~~\left[ e_{3},e_{1}\right] =0
\end{equation*}%
Note that $M^{3}$ is not flat since the scalar curvature is $r=-10-\frac{%
1+8x^{2}}{2\left( x^{2}\right) ^{2}}$.

\begin{remark}
Given an associated metric $g$, any Riemann metric $\gamma $ can be obtained
by adding to $g$ a tensor field $t$, provided that $\gamma $ remains a
symmetric, positive definite bilinear form. The following formula gives such
a $t$ in simplifying coordinates:%
\begin{equation}
t=\left[ 
\begin{array}{ccc}
e^{\iota }-1 & \rho  & \sigma  \\ 
\rho  & e^{\kappa }-\frac{1+\alpha ^{2}}{\beta ^{2}} & \upsilon  \\ 
\sigma  & \upsilon  & e^{\nu }-\frac{\beta ^{2}\left( 1+F^{2}\right)
+\epsilon ^{2}}{\beta ^{2}\zeta ^{2}}%
\end{array}%
\right]   \label{mt}
\end{equation}%
where $\iota $, $\kappa $, $\nu $, $\rho $, $\sigma $ and $\upsilon $ are
any smooth functions of $x^{1}$, $x^{2}$ and $x^{3}$. The final formula
giving any Riemannian metric is%
\begin{equation}
g+t=\gamma =\left[ 
\begin{array}{ccc}
e^{\iota } & -\frac{\alpha }{\beta }+\rho  & -\frac{F}{\zeta }+\sigma  \\ 
-\frac{\alpha }{\beta }+\rho  & e^{\kappa } & \frac{\alpha \beta F-\epsilon 
}{\beta ^{2}\zeta }+\upsilon  \\ 
-\frac{F}{\zeta }+\sigma  & \frac{\alpha \beta F-\epsilon }{\beta ^{2}\zeta }%
+\upsilon  & e^{\nu }%
\end{array}%
\right]   \label{rm}
\end{equation}%
Recalling step (\ref{step2}) of the algorithm presented in section \ref%
{algorithm}, we see that in (\ref{rm}) the functions $\alpha $, $\beta $, $%
\zeta $, $\epsilon $ and $F$ can be absorbed by a redefinition of $\rho $, $%
\sigma $ and $\upsilon $ thus yielding (as expected) any symmetric, positive
definite bilinear form as a Riemannian metric. However in (\ref{mt}), the
matrix $t$ involves the functions $\alpha $, $\beta $, $\zeta $, $\epsilon $
and $F$ in a crucial way.
\end{remark}

\end{document}